\documentclass[aap,MSNbibl,nameyear,seceqn,dvips]{arximspdf}

\doi{10.1214/12-AAP882} 
\volume{23}
\issue{4}
\pubyear{2013}
\firstpage{1629}
\lastpage{1659}

\makeatletter
\newcommand{\eqref}[1]{(\ref{#1})}
\def\one{\mathbf{1}}
\def\reff#1{(\ref{#1})}
\def\P{\mathbb{P}}
\def\R{\mathbb{R}}
\def\Z{\mathbb{Z}}
\def\N{\mathbb{N}}
\def\LL{\mathcal{G}}

\newtheorem{theo}{Theorem}
\newtheorem{prop}{Proposition}
\newtheorem{lem}{Lemma}
\newproclaim{defin}{Definition}
\newtheorem{cor}{Corollary}
\newproclaim{ex}{Example}
\newproclaim{algo}{Algorithm}
\makeatother

\begin{document}
\begin{frontmatter}

\title{Kalikow-type decomposition for multicolor infinite range particle systems}
\runtitle{Multicolor infinite range particle systems}

\begin{aug}
\author[A]{\fnms{A.} \snm{Galves}\ead[label=e1]{galves@usp.br}\thanksref{t1}},
\author[B]{\fnms{N.~L.} \snm{Garcia}\corref{}\ead[label=e2]{nancy@ime.unicamp.br}\thanksref{t2}},
\author[C]{\fnms{E.} \snm{L\"ocherbach}\ead[label=e3]{eva.loecherbach@u-cergy.fr}\thanksref{t3}}
\and
\author[D]{\fnms{E.} \snm{Orlandi}\ead[label=e4]{orlandi@mat.uniroma3.it}\thanksref{t4}}
\thankstext{t1}{Supported in part by a CNPq fellowship Grant 305447/2008-4.}
\thankstext{t2}{Supported in part by a CNPq fellowship Grant 302755/2010-1.}
\thankstext{t3}{Supported by ANR-08-BLAN-0220-01 and Galileo 2008-09.}
\thankstext{t4}{Supported by Prin07: 20078XYHYS.}
\runauthor{A. Galves et al.}
\affiliation{IME-USP, IMECC/UNICAMP,
Universit\'e de Cergy-Pontoise and~Universit\`a~di~Roma~Tre}
\address[A]{A. Galves\\
IME-USP\\
Rua do Mat\~{a}o, 1010---Cidade Universit\'{a}ria\\
05508-090---S\~{a}o Paulo---SP\\
Brazil\\
\printead{e1}} 
\address[B]{N.~L. Garcia\\
IMECC/UNICAMP\\
Rua S\'{e}rgio Buarque de Holanda\\
651---Cidade Universit\'{a}ria\\
13.083-859---Campinas---SP\\
Brazil\\
\printead{e2}}
\address[C]{E. L\"ocherbach\\
Laboratoire AGM\\
CNRS UMR 8088\\
Universit\'{e} de Cergy-Pontoise\\
2, avenue Adolphe Chauvin\\
95302 Cergy-Pontoise Cedex\\
France\\
\printead{e3}}
\address[D]{E. Orlandi\\
Dipartimento di Matematica\\
Universit\`{a} di Roma Tre\\
Largo S. Leonardo Murialdo 1\hspace*{5.2pt}\\
I-00146 Roma\\
Italy\\
\printead{e4}}
\end{aug}

\received{\smonth{8} \syear{2010}}
\revised{\smonth{5} \syear{2012}}

%
\begin{abstract}
We consider a particle system on $\Z^d$ with real state space and
interactions of infinite range. Assuming that the rate of change is
continuous we obtain a Kalikow-type decomposition of the infinite
range change rates as a mixture of finite range change rates.
Furthermore, if a high noise condition holds, as an application of
this decomposition, we design a feasible perfect simulation
algorithm to sample from the stationary process. Finally, the
perfect simulation scheme allows us to forge an algorithm to obtain an
explicit construction of a coupling attaining Ornstein's $\bar
d$-distance for two ordered Ising probability measures.
\end{abstract}

%

\end{frontmatter}

\section{Introduction}

In this paper we present a Kalikow-type decomposition for interacting
multicolor systems on $\Z^d$ having real state space and
interactions of infinite range. By a Kalikow-type decomposition we
mean a representation of the infinite range rates as a countable
mixture of local change rates of increasing range. This decomposition
extends the notion of random Markov chains to interacting particle
systems and has many potential theoretical consequences and
applications. As a first example we present a perfect simulation
algorithm which is based on the decomposition. As a corollary we obtain a
result about the existence and uniqueness of the invariant measure of
the system as well as a rate of convergence to stationarity. As a second
application we construct a coupling attaining the
$\bar{d}$-distance for two ordered Ising probability measures.

By a perfect simulation algorithm we mean a simulation which samples
in a finite window precisely from the stationary law of the infinite
process. More precisely, for any finite set of sites $F$ we want to
sample the projection of the stationary law on $F$. Our approach is
feasible in the sense that it stops almost surely after a finite
number of steps. It does not require any duality or monotonicity properties.
We do not assume that the system has a dual, or is attractive, or
monotone in any sense. Our system is not spatially homogeneous. The
basic assumptions are the continuity of the infinite range change
rates together with a high-noise condition
[Condition (\ref{eqcondition2}): fast decay of the range
influence on the change rate and a certain subcriticality-criterion].

Concerning possible applications, perfect simulation of infinite range
continuous (or discrete) systems has shown to be an important tool,
for instance, in statistical inference for Gibbs distributions, for
Bayesian statistics, for maximum likelihood estimation, rates of
convergence of estimators. This field of research has enormous
relevance due to its applications in image processing, spatial
statistics, gene expression, to cite just a few of the possible
applications. Some of the applications use discrete state space
whereas others deal with continuous state space, some deal with
nearest neighbor interaction (image processing) whereas others need
infinite range interactions (tumor growth, dynamics of populations of
neurons or gene expression). We refer the reader to the books
of \citet{MllWaa04}, \citet{GaeGuy10} for
some applications.

Let us stress that from an applied point of view, it is important to
be able to deal with infinite range systems. In neuroscience, for
example, for all practical purposes, the interactions between neurons
have infinite range [see, e.g., \citet{Cesetal09} and \citet{CesNasVas10}]. Also, very
simple error structures superposed on finite range systems will
produce infinite range systems. This is the case, for example, of
blurred images [see, e.g., \citet{NisWon99} and \citet{Tan02}].
Up to now, denoising of images has been studied mostly with
very simple models (finite lattice, nearest neighbor interaction) due
to the difficulty to sample from a Gibbs distribution with infinite
range interactions and/or continuous state space.

So it is important to build up a unifying approach that enables us to
deal with infinite range interactions and/or continuous state space at
the same time. Although perfect simulation of Gibbs random fields are
known for discrete state spaces and finite range interaction, very
little is known for simulation of continuous state spaces and/or
infinite range interaction. To the best of our knowledge
ours is the first concrete result concerning perfect simulation for
interacting particle systems with continuous state space and infinite
range interactions.

This paper is organized as follows. The model and the Kalikow-type
decomposition (Theorem~\ref{theodecomp}) are presented in Section
\ref{secdecomp}. The aim of Section~\ref{ex} is to present some examples where the
decomposition can be explicitly done. In particular, we apply
Theorem~\ref{theodecomp} to the important case of Gibbs measures with
infinite range interactions and continuous spin values.
In Section~\ref{perfect} we present the perfect simulation algorithm as a main
application of the Kalikow-type decomposition. In particular, Theorem
\ref{theo2} shows that the proposed algorithm is feasible under a high
noise condition. The proofs are given in Sections~\ref{secproof1}
and~\ref{secbw}. In Section~\ref{sectionappli} we give as a main
application of the perfect simulation algorithm an explicit
construction of a coupling attaining Ornstein's $\bar d$-distance for
two ordered Ising probability measures. We present a section
discussing the user impatient bias. We conclude the article with
final comments and a bibliographical discussion.

\section{Definitions, notation and convex
decomposition} \label{secdecomp}

We consider interacting particle systems on $\Z^d$ having
state space $A$ and interactions of infinite range. The elements of the
state space $A$ are called \textit{colors}. To each site in $\Z^d$ we assign
a color. The coloring of the sites changes as time goes by. The rate
at which the color of a fixed site $i$ changes from a color $a$ to a
new color $b$ is a function of the entire configuration and depends on
$b$.

In what follows, we suppose that $ A$ is a Borel subset of $\R, $
equipped with its
Borel sigma-field $\mathcal{A}$.
$\varrho$ will be a finite nonnegative reference measure on $(A,
\mathcal{A} ). $
The initial lowercase letters $a$, $b$, $c, \ldots$ will denote
elements of~$A.$
We denote by $S= A^{\Z^d}$ the configuration space with its product sigma
algebra,~$\mathcal{S} $. We call an element of $S $ a configuration.
Configurations will be denoted by Greek letters
$\eta, \zeta, \xi, \ldots.$ A point $ i \in\Z^d $ will be
called a site. We define on $\Z^d$ the $L^1$ norm, $\| i \|= \sum_{k=1}^d |i_k|$.
For $k \ge0$, let the ball of radius $k$ be
\[
V_{i} (k) = \bigl\{j \in\Z^d; \|j - i \| \le k\bigr\}.
\]

As usual, for any $i \in\Z^d$, $\eta(i)$ will denote the
value of the configuration $\eta$ at site~$i$. By extension, for any
subset $V \subset\Z^d$, $\eta(V)\in A^V$ will denote the restriction
of the configuration $\eta$ to the set of positions in $V.$ For any
$\eta,$ $i$ and $a,$ we shall denote $\eta^{i,a}$ the modified
configuration
\[
\eta^{i,a}(j) = \eta(j) \qquad\mbox{for all $j \neq i,$ and $
\eta^{i,a}(i) = a.$}
\]
For any $i \in\Z^d,$ let $ c_i (a , \eta) $ be a
positive $\mathcal{A}\otimes\mathcal{S} - \mathcal{B} ( \R_+)$-measurable function such that
the following two properties hold. First, for $\varrho$-almost all $ a
\in A, $ $ \eta\mapsto c_i (a , \eta) $
is continuous. Second, we have
%
\begin{equation}
\label{eqboundedrate} \sup_{ i \in\Z^d } \sup_\eta\int
_A c_i (a, \eta) \varrho(da) < \infty.
\end{equation}

A multicolor system with interactions of infinite range is a Markov
process on~$S $ whose generator is defined on cylinder functions by
%
\begin{equation}
\label{eqgenerator} \LL f(\eta) = \sum_{i \in\Z^d} \int
_{A} \varrho(da) c_i ( a, \eta) \bigl[f\bigl(
\eta^{i,a}\bigr) - f(\eta)\bigr] ,
\end{equation}
where $f \in D ( \LL) = \{ f \dvtx  |\!|\!| f |\!|\!| = \sum_{i \in\Z^d } \Delta_f (i) < \infty\} $ with
$ \Delta_f (i ) = \sup\{ | f(\eta) - f (\zeta)| \dvtx  \eta(j) = \zeta
(j) \mbox{ for all } j \neq i \} .$

By Theorem 3.9 of Chapter 1 of \citet{Lig85N1} the following condition,
together with (\ref{eqboundedrate}),
implies that $\LL$ is the generator of a Feller process $(\sigma_t )$
on $S$:
%
\begin{equation}
\label{A3} \sup_{i \in\Z^d} \sum_{j \neq i }
\sup_\eta \sup_{b \in A} \biggl\{ \int_A
\rho(da) \bigl| c_i (a,\eta)- c_i \bigl(a,\eta^{j, b}
\bigr)\bigr| \biggr\} < \infty.
\end{equation}
In the following we shall work under conditions stronger than
(\ref{A3}) ensuring not only that $ \LL$ is the generator of a unique
Feller process, but also the possibility of perfectly simulating the
stationary distribution corresponding to this infinitesimal generator.
As a
byproduct this implies that the system admits the existence of a
unique invariant measure $\mu$.

The main result of this article is a Kalikow-type convex decomposition
of the change rates. We will prove that the change rate can be
decomposed as
%
\begin{equation}
\label{eqconvex1} c_i (a , \eta) = M_i \biggl[
\lambda_i(-1) p^{[-1]}_i(a) + \sum
_{ k
\geq0} \lambda_i ( k) p_i^{[k]}
\bigl( a | \eta\bigl(V_i (k)\bigr) \bigr) \biggr],
\end{equation}
where:
\begin{itemize}
\item$M_i, i \in\Z^d$ are positive constants,
\item for each $i \in\Z^d$, $\{\lambda_i(k), k \ge-1\}$ is a
probability distribution on $\{-1, 0, 1,\break 2, \ldots\}$,
\item for each $i \in\Z^d$, $p_i^{[-1]}(\cdot)$ is a probability
density on $A$ with respect to the reference measure $\varrho,$ which does
not depend on the configuration,
\item for each $k \ge0$ and for each $\eta\in S ,$ $p_i^{[k]}(
\cdot| \eta( V_i (k)))$ is a probability density with respect to the
reference measure $\varrho,$ depending only on the local
configuration $ \eta( V_i (k))$.
\end{itemize}

For convenience of the presentation we will add additional invisible
jumps in~\eqref{eqgenerator}. This is obtained by adding a cemetery
$\Delta$ to $A $ and defining $A^* := A \cup\{ \Delta\} .$ Define also
\[
\varrho^* := \varrho+ \delta_{\Delta}.
\]

Denote
%
\begin{equation}
\label{eqidiot} M_i := \sup_{\eta\in A^{\Z^d}} \int c_i (a,
\eta ) \varrho(da).
\end{equation}
Notice that $M_i$ is finite under condition (\ref{eqboundedrate}),
and define
%
\begin{equation}
\label{eqdelta} c_i ( \Delta, \eta) := M_i - \int
_A c_i ( a, \eta) \varrho(da) .
\end{equation}
Observe that
%
\begin{equation}
\label{eqnull} \inf_{\eta} c_i (\Delta, \eta) = 0.
\end{equation}
Therefore we can rewrite the generator given by \eqref{eqgenerator} as
%
\begin{equation}
\label{eqgeneratorbis} \LL f(\eta) = \sum_{i \in\Z^d}
\int_{A^* } \varrho^* (da) c_i ( a, \eta) \bigl[f
\bigl(\eta^{i,a}\bigr) - f(\eta)\bigr] ,
\end{equation}
where, by convention, for any $i \in\Z^d$ and any $\eta\in S=A^{\Z
^d}$, we define
\[
\eta^{i, \Delta} = \eta.
\]
It follows that \eqref{eqgeneratorbis} is a representation of the
same generator as \eqref{eqgenerator}.

In order to obtain the
decomposition we need the following continuity condition.

\textit{Continuity condition}.
%
\begin{equation}
\label{eqcontinuity} \sup_{i \in\Z^d} \int_A
\sup_{\eta, \zeta: \eta(V_i(k)) =
\zeta(V_i(k))} \bigl| c_i (a, \eta) - c_i (a, \zeta)\bigr|
\varrho(da) \rightarrow0
\end{equation}
as $k \rightarrow\infty.$

To describe the convex decomposition of the rate function $c_i$, we
have to introduce the following quantities. Define
%
\begin{equation}
\label{eqalpha0} \alpha_{i} (-1) = \int_{A^* }
\inf_{\zeta\in A^{\Z^d}} c_i(a, \zeta) \varrho^* (da) ,
\end{equation}
and for any $k \ge0 , $
%
\begin{equation}
\label{eqalpha} \alpha_i (k) = \inf_{w \in A^{V_i (k) }} \biggl( \int
_{A^* } \inf_{\zeta: \zeta(V_i (k)) = w} c_i( a,\zeta)
\varrho^* (da ) \biggr) .
\end{equation}
The continuity of $ c_i (a, \eta)$ in $\eta $ and the separability of
$S$ imply the measurability of
$\inf_{\zeta: \zeta(V_i (k)) = w} c_i( a,\zeta) $ and $\inf_{\zeta
\in A^{\Z^d}} c_i(a, \zeta)$ with respect to $a.$

Note that by \eqref{eqnull}
\[
\int_{A^* } \inf_{\zeta\in A^{\Z^d}} c_i(a, \zeta)
\varrho^* (da)= \int_{A } \inf_{\zeta\in A^{\Z^d}}
c_i(a, \zeta) \varrho (da).
\]
Further,
by construction, we have that $\alpha_i (k) \le\alpha_i (k + 1), $
for each $k \geq- 1$. We claim that
%
\begin{equation}
\label{eqm} M_i = \lim_{k \to\infty} \alpha_i (k).
\end{equation}

To obtain equality (\ref{eqm}), fix some $w \in A^{V_i (k) }$;
from (\ref{eqdelta}), we have that
\begin{eqnarray*}
&&\int_{A^* } \inf_{\zeta: \zeta(V_i (k)) = w} c_i( a,
\zeta) \varrho^* (da )
\\
&&\qquad = \int_{A} \inf_{\zeta: \zeta(V_i (k)) = w} c_i( a,
\zeta) \varrho^* (da ) + \inf_{\zeta: \zeta(V_i (k)) = w} c_i(\Delta,\zeta)
\\
&&\qquad = \int_A \inf_{\zeta: \zeta(V_i (k)) = w} c_i( a,
\zeta) \varrho(da) + M_i - \sup_{\zeta: \zeta(V_i (k)) = w} \int
_A c_i (a, \zeta) \varrho(da) .
\end{eqnarray*}
But
\[
\int_A \inf_{\zeta: \zeta(V_i (k)) = w} c_i( a,\zeta)
\varrho(da) -\sup_{\zeta: \zeta(V_i (k)) = w} \int_A c_i
(a, \zeta) \varrho (da) \to0
\]
as $k \to\infty$ thanks to condition \eqref{eqcontinuity}.

Hence, to each
site $i$ we can associate a probability distribution $ \lambda_i$ by
%
\begin{equation}
\label{eqlambda0} \lambda_{i} (-1) = \frac{\alpha_{i} (-1)}{M_i} ,
\end{equation}
and for $k \ge0 $
%
\begin{equation}
\label{eqlambdak} \lambda_i (k) = \frac{\alpha_i (k)  -  \alpha_{i} (k-1)}{M_i } .
\end{equation}

Now we are ready to state the decomposition theorem.

%
\begin{theo} \label{theodecomp} Let $(c_i)_{i \in\Z^d} $ be a family
of measurable rate functions satisfying conditions
$(\ref{eqboundedrate}), (\ref{eqdelta})$ and $(\ref{eqcontinuity})$.
Then, for each site $i$, for $M_i$ defined by $(\ref{eqidiot})$
and $\lambda_i(\cdot)$ defined by $(\ref{eqlambda0})$ and
$(\ref{eqlambdak})$, there exist:
\begin{itemize}
\item$p_i^{[-1]}$ a probability density with respect to $\varrho$
with support $A $,
\item a family of conditional
probability densities $p_i^{[k]}$ $[$given by $(\ref{eqpik})]$, $k
\ge0$ on $A^*,$ with respect to
$\varrho^*$, depending on the local configurations $\eta(V_i (k))
\in A^{V_i (k)}$
such that
%
\begin{equation}
\label{cmmc1} c_i (a, \eta) = M_i p_i (a|
\eta)\qquad \mbox{for $\varrho^*$-almost all $a \in A^* ,$}
\end{equation}
where
%
\begin{equation}
\label{cmmc2} p_i(a|\eta) = \lambda_i(-1)
p^{[-1]}_i(a) + \sum_{ k
\geq0}
\lambda_i ( k) p_i^{[k]} \bigl( a | \eta
\bigl(V_i (k)\bigr) \bigr).
\end{equation}
\end{itemize}

As a consequence, the infinitesimal generator $\LL$ given by
$(\ref{eqgeneratorbis})$ can be rewritten as
%
\begin{eqnarray}\quad
\label{eqgenerator2} \LL f(\eta) &=& \sum_{i \in\Z^d}
M_i \biggl[ \lambda_i (- 1) \int_A
p_i^{[-1]} (a) \bigl[f\bigl(\eta^{i,a}\bigr) - f(
\eta)\bigr] \varrho(da)
\nonumber
\\[-8pt]
\\[-8pt]
\nonumber
&&\hspace*{38pt} {} + \sum_{ k \geq0 } \lambda_i (k)
\int_{A^*} p_i^{[k]}\bigl(a| \eta
\bigl(V_i (k)\bigr)\bigr) \bigl[f\bigl(\eta^{i,a}\bigr) - f(
\eta)\bigr] \varrho^* (da) \biggr] .
\nonumber
\end{eqnarray}
\end{theo}

Note that for $k = -1,$ $p_i^{[-1]} (a)$ does not depend on the
configuration and $\lambda_{i} (-1)$ represents the spontaneous
self-coloring rate of site $i$ in the process. We will see in the
proof that $p_i^{[-1]}$ is defined in such way that $p_i^{[-1]}
(\Delta) = 0$ and therefore, the choice $k = - 1 $ always implies a
choice of a real color $ a \in A , $ not of $ a = \Delta.$

The decomposition given in Theorem~\ref{theodecomp} was designed in
such way that the probability of self-coloring is maximized. This is
important to speed up the perfect simulation algorithm. Obviously,
slight modifications can be employed for different purposes as we will
see in Example~\ref{exam2} (Section~\ref{ex}).

The representation given by (\ref{eqgenerator2}) provides a random
finite range description of the time evolution of the process.
We start with an initial configuration $\eta$ at time zero. For each
site $ i \in\Z^d ,$ we consider a rate $M_i$ Poisson point process
$N^i .$ The Poisson processes corresponding to distinct sites are all
independent. If at time $t$, the Poisson clock associated to site $i$
rings, we choose a range $k$ with probability $\lambda_i (k)$
independently of everything else. Then, we update the value of
the configuration at this site by choosing a symbol $a$ with probability
$p_i^{[k]} (a | \sigma_t(V_i(k))) \varrho^* (da) $. Choosing
the symbol $\Delta$ means that we actually keep the current value of
the spin.

In Section~\ref{ex} we give examples of infinite range interacting
systems where Theorem~\ref{theodecomp} can be applied.

\section{\texorpdfstring{Proof of Theorem \protect\ref{theodecomp}}{Proof of Theorem 1}} \label{secproof1}

Put for any $a \in A^* ,$
\begin{eqnarray*}
c_i^{[-1]} (a) &=& \inf_{\zeta} c_i(a ,
\zeta),
\\
\Delta^{[-1]}_i (a) &= &c_i^{[-1]} (a),
\\
c_i^{[0]} \bigl( a| \eta(i)\bigr)& =& \inf_{\zeta: \zeta(i ) =
\eta(i ) }
c_i (a, \zeta) ,
\\
\Delta^{[0]}_i \bigl( a | \eta(i )\bigr)&= &
c_i^{[0]} \bigl( a |\eta(i )\bigr) - c^{[-1]}_i
(a) .
\end{eqnarray*}
For any $k \geq1,$ define
\begin{eqnarray*}
c_i^{[k]} \bigl( a| \eta\bigl(V_i(k)\bigr)
\bigr) &= &\inf_{ \zeta: \zeta(V_i (k)) =
\eta(V_i (k))} c_i(a , \zeta) ,
\\
\Delta_i^{[k]} \bigl(a | \eta\bigl(V_i (k)
\bigr)\bigr)&=& c_i^{[k]} \bigl(a |\eta\bigl(V_i
(k)\bigr)\bigr) - c^{[k-1]}_{i} \bigl( a | \eta\bigl(
V_i ({k-1}) \bigr)\bigr).
\end{eqnarray*}

Then we have that for any $ a \in A ,$
%
\begin{equation}
\label{eqcia} c_i(a, \eta) = \sum_{j=-1}^k
\Delta^{[j]}_i \bigl(a| \eta\bigl(V_i(j)\bigr)
\bigr) + \bigl[ c_i(a,\eta) - c_i^{[k]} \bigl(a|
\eta\bigl(V_i(k)\bigr)\bigr) \bigr].
\end{equation}

Note that
\[
c^{[-1]}_i (\Delta) = \inf_\eta c_i (
\Delta, \eta) = M_i - \sup_\eta\int_A
c_i (a, \eta) \varrho(da ) = 0.
\]
Therefore, for $a=\Delta$ decomposition (\ref{eqcia}) starts with $j
= 0$,
\[
c_i ( \Delta, \eta) = \sum_{j= 0 }^k
\Delta^{[j]}_i \bigl(\Delta| \eta\bigl(V_i(j)
\bigr)\bigr) + \bigl[ c_i(\Delta,\eta) - c_i^{[k]}
\bigl(\Delta|\eta\bigl(V_i(k)\bigr)\bigr) \bigr].
\]

By monotonicity, we have for $\varrho^*$-almost all $a \in A^*$
that
\[
c_i^{[k]} \bigl(a|\eta\bigl(V_i(k)\bigr)\bigr)
\to \lim_{k } c_i^{[k]} \bigl(a|\eta
\bigl(V_i(k)\bigr)\bigr) \le c_i ( a, \eta) \qquad\mbox{as } k
\to\infty.
\]
Hence, by monotone convergence,
\[
\int_A c_i^{[k]} \bigl(a|\eta
\bigl(V_i(k)\bigr)\bigr) \varrho(da ) \to \int_A
\Bigl[ \lim_{k } c_i^{[k]} \bigl(a|\eta
\bigl(V_i(k)\bigr)\bigr) \Bigr] \varrho(da) \le\int
_A c_i (a, \eta) \varrho(da) .
\]
On the other hand, by (\ref{eqcontinuity}),
\[
\int_A c_i^{[k]} \bigl(a|\eta
\bigl(V_i(k)\bigr)\bigr) \varrho(da ) \to\int_A
c_i (a, \eta) \varrho(da) .
\]
Hence, for $\varrho^*$-almost all $a$ and all $ \eta,$
\[
\lim_{k } c_i^{[k]} \bigl(a|\eta
\bigl(V_i(k)\bigr)\bigr) = \sum_{j=-1}^\infty
\Delta_i^{[j]} \bigl(a| \eta\bigl(V_i (j)\bigr)
\bigr) = c_i(a, \eta) .
\]
Taking into account \eqref{eqalpha0} and \eqref{eqlambda0},
\[
M_i \lambda_i ( - 1) = \int_A
\Delta_i^{[-1]} ( a) \varrho(da ).
\]
Hence, we can define
\[
p_i^{[-1]} (a) = \frac{\Delta_i^{[-1]} (a)}{ M_i \lambda_i (-1)}
\]
and
\[
p_i^{[-1]} ( \Delta) = 0 .
\]
Hence, $p_i^{[-1]} (a) $ is a probability density with respect to
$\varrho^* .$
Now, for $ k \geq0,$ put
%
\begin{equation}
\label{eqlambdaalmost} \tilde{\lambda}_i \bigl(k,\eta
\bigl(V_i(k)\bigr)\bigr) = \frac{1}{M_i} \int_{A^*}
\Delta_i^{[k]} \bigl(a| \eta\bigl(V_i (k)\bigr)
\bigr) \varrho^* (da) ,
\end{equation}
and for any $i,k$ such that $ \tilde{\lambda}_i (k,\eta(V_i(k))) >
0,$ we define
\[
\tilde{p}_i^{[k]} \bigl(a | \eta\bigl(V_i(k)
\bigr)\bigr) = \frac{\Delta_i^{[k]} (a|
\eta(V_i (k)))}{ M_i   \tilde{\lambda}_i (k,\eta(V_i(k)))}.
\]
For $i, k$ such that
$ \tilde{\lambda}_i (k,\eta(V_i(k))) = 0,$ define $\tilde
{p}_i^{[k]} (a |
\eta(V_i(k)))$ in an arbitrary fixed way.\vadjust{\goodbreak}

Hence, for $\varrho^*$-almost all $a \in A^*,$
%
\begin{equation}
\label{eqalmost}\qquad c_i(a,\eta) = M_i \Biggl[
\lambda_i (- 1) p_i^{[-1]} (a) + \sum
_{k=
0 }^\infty\tilde{\lambda}_i \bigl(k,\eta
\bigl(V_i(k)\bigr)\bigr) \tilde{p}_i^{[k]}
\bigl(a| \eta\bigl(V_i(k)\bigr)\bigr) \Biggr] .
\end{equation}
In (\ref{eqalmost}) the factors $\tilde{\lambda}_i (k,\eta(V_i(k) )),
k \geq0, $ still depend on $\eta(V_i(k)) .$ To obtain the
decomposition as in the theorem, we must rewrite it as follows.

For any $i,$ take $M_i$ as in (\ref{eqm}) and the sequences $\alpha_i
(k), \lambda_i (k), k \geq-1,$ as defined in (\ref{eqalpha}) and
(\ref{eqlambdak}), respectively. Define the new quantities
\[
\alpha_i \bigl(k,\eta\bigl(V_i(k)\bigr)\bigr) =
M_i \sum_{l \le k} \tilde{\lambda
}_i \bigl(l, \eta\bigl( V_i(l)\bigr)\bigr)
\]
and notice that
\[
\alpha_i \bigl( k,\eta\bigl(V_i(k)\bigr)\bigr) = \int
_{A^*} c_i^{[k]} \bigl(a, \eta
\bigl(V_i(k)\bigr)\bigr) \varrho^* (da)
\]
is the total mass associated to $ c_i^{[k]} ( \cdot,\eta(V_i(k))) .$

By definition of $ \alpha_i (k) $
in (\ref{eqalpha}), $ \alpha_i (k) $ is the smallest total mass
associated to $ c_i^{[k]},$
uniformly with respect to all possible neighborhoods $ \eta( V_i(k)).$
Hence, in order to get
a decomposition with weights $ \lambda_i (k) $ not depending on the
configuration, we have to define a partition of the interval $ [ 0,
\alpha_i (k, \eta(V_i(k)))]$
according to the values of $ \alpha_i (k) .$

This yields, for any $k \geq0$, the following definition of the
conditional finite range
probability densities.
%
\begin{eqnarray}
\label{eqpik}
&& p_i^{[k]} \bigl( a| \eta
\bigl(V_i({k})\bigr)\bigr)
\nonumber
\\[-2pt]
&&\qquad= \sum_{-1 = l' \le l }^{k-1} 1_{\{
\alpha_i (l' - 1 ,\eta(V_i (l' -1))) < \alpha_i ({k-1}) \le\alpha_i({l'},
\eta(V_i({l' })))\}}\nonumber\\[-2pt]
&&\qquad\hspace*{45pt}\times 1_{\{\alpha_i (l,\eta(V_i (l))) < \alpha_i ({k}) \le\alpha_i({l+1}
\eta(V_i({l+1})))\}}
\nonumber
\\[-8pt]
\\[-8pt]
\nonumber
&&\hspace*{45pt}\qquad{}\times \Biggl[ \frac{\alpha_i (l',\eta(V_i(l'))) - \alpha_i(k-1) }{M_i
\lambda_i({k})} \tilde{p}_i^{[l']} \bigl(a |
\eta\bigl(V_i\bigl(l'\bigr)\bigr)\bigr)
\\[-2pt]
&&\hspace*{50pt}\quad\qquad{} + \sum_{m = l'+1}^{l} \frac{\tilde{\lambda}_i
(m , \eta(V_i (m)))}{ M_i \lambda_i (k) }
\tilde{p}_i^{[m]} \bigl(a | \eta\bigl(V_i(m)
\bigr)\bigr)
\nonumber
\\[-2pt]
&&\hspace*{50pt}\quad\qquad{} + \frac{\alpha_i({k}) - \alpha_i
(l,\eta(V_i(l)))}{M_i
\lambda_i({k})} \tilde{p}_i^{[l+1]} \bigl(a|
\eta\bigl(V_i({l+1})\bigr)\bigr) \Biggr].\nonumber
\end{eqnarray}

The desired decomposition now follows from this.

\section{Examples}\label{ex}

In this section we give examples where the decomposition of Theorem
\ref{theodecomp} can be applied. We start with an example from
Bayesian statistics and image reconstruction.\vadjust{\goodbreak}

\begin{ex}[(Autonormal distribution)] This model can be seen
as the spatial analogue of the autoregressive model. The usual way to
describe its dynamics is through the simultaneous schemes: Each pixel
updates its value using a normal distribution with mean depending on
the values of its neighbors.
\end{ex}

In this work, we are going to generalize this definition to
incorporate long-range interactions and arbitrary neighborhoods, but
we are going to limit the values of the process to be on a compact
interval $[l,u]$. The existence of a unique invariant measure for this
process was studied by \citet{McBSpe77} and revisited by
\citet{FerGry08}. \citet{Gib04}
proposes a finite version of this model as the posterior distribution
for Bayesian restoration of grayscale images. \citet{Hub07} studies
perfect simulation of these distributions in a finite box.

Let $\{\sigma(i), i \in\Z^d\}$ be a collection of positive real
numbers. Consider that each pixel in $\Z^d$ has an independent
exponential clock and when the clock rings at the pixel $i$ it updates
its value depending on the values of its neighbors using a normal
distribution with mean $h(i,\eta)$ and variance $\sigma^2(i)$
conditioned to lie in a given compact interval. Without loss
of generality we can consider the interval $[0,1]$. The term
$h(i,\eta)$ depends on $\eta$ only through the values on the
neighborhood of the site $i$. Typically,
\[
h(i,\eta) = \sum_{j \neq i} J(i-j) \eta(j),
\]
where $J\dvtx  \Z^d \rightarrow\R^+$ is summable, nonnegative and
symmetric: $J(i) \ge0$ for all $i \in\Z^d$, $J(0) = 0$ and $0 < J :=
\sum_{i \in\Z^d} J(i) =1$.

In this case, $\varrho( da) = \one_{[0,1]} (a) \,da $ and
%
\begin{eqnarray}
\label{eqanci} c_i(a, \eta) = \frac{1}{\sigma(i)} \frac{\phi((a-h(i,\eta))/\sigma(i))}{\Phi((1-h(i,\eta))/\sigma
(i)) - \Phi(-h(i,\eta)/\sigma(i))},
\nonumber
\\[-8pt]
\\[-8pt]
\eqntext{0 \le a \le1,}
\end{eqnarray}
where
$\phi$ and $\Phi$ are the density and cumulative function of the
standard normal distribution, respectively.

Applying the bounds given by Proposition 1 of Fern\'andez, Ferrari and
Grynberg (\citeyear{FerFerGry07}) we can show that (\ref{eqanci}) satisfies the
assumptions (\ref{eqboundedrate}) and (\ref{eqcontinuity}) needed in
order to apply the decomposition of Theorem~\ref{theodecomp}. In our
case,
$M_i = 1$ and
%
\begin{eqnarray}
\label{eqanalpha1}
\alpha_i (-1)& =& \frac{1}{A^{+}} \biggl(\Phi
\biggl(\frac{x_i(\sigma
(i)) -\mu^{+}}{\sigma(i)} \biggr) - \Phi \biggl(\frac{-\mu^{+}}{\sigma(i)} \biggr) \biggr)
\nonumber
\\[-8pt]
\\[-8pt]
\nonumber
&&{}+\frac{1}{A^{-}} \biggl(\Phi \biggl(\frac{1 -
\mu^{-}}{\sigma(i)} \biggr) - \Phi \biggl(
\frac{x_i(\sigma(i))-\mu^{-}}{\sigma(i)} \biggr) \biggr),
\end{eqnarray}
where $\mu^{-} = \inf_{\eta} h(i,\eta)$, $\mu^{+} = \sup_{\eta}
h(i,\eta)$,
\[
A^{-} = \Phi \biggl( \frac{1 - \mu^{-}}{\sigma(i)} \biggr) - \Phi \biggl(
\frac{ - \mu^{-}}{\sigma(i)} \biggr),\qquad A^{+} = \Phi \biggl( \frac{1 - \mu^{+}}{\sigma(i)}
\biggr) - \Phi \biggl( \frac{ - \mu^{+}}{\sqrt{\sigma(i)}} \biggr)
\]
and
\[
x_i\bigl(\sigma(i)\bigr) = \frac{\mu^{-} + \mu^{+}}{2} - \frac{\sigma
(i)^2}{\mu^{+}-
\mu^{-}}
\log\frac{A^{-}}{A^{+}}.
\]

Also,
%
\begin{eqnarray}
\label{eqanalpha1}
\alpha_i (k) &=& \inf_{w \in A^{V_i (k) }}
\frac{1}{A^{+}_{i}(k)} \biggl(\Phi \biggl(\frac{x_i(k,\sigma(i)) -
\mu^{+}_{i}(k)}{\sigma(i)} \biggr) - \Phi \biggl(
\frac{-\mu^{+}_{i}(k)}{\sigma(i)} \biggr) \biggr)
\nonumber
\\[-8pt]
\\[-8pt]
\nonumber
&&\hspace*{34pt}{}+\frac{1}{A^{-}_{i}(k)} \biggl(\Phi \biggl(\frac{1 -
\mu^{-}_{i}(k)}{\sigma(i)} \biggr) - \Phi \biggl(
\frac{x_i(k,\sigma(i))-\mu^{-}_{i}(k)}{\sigma(i)} \biggr) \biggr),
\end{eqnarray}
where $\mu^{-}_{i}(k) = \inf_{\eta: \eta(V_i (k)) = w} h(i,\eta)$,
$\mu^{+}_{i}(k) = \sup_{\eta: \eta(V_i (k)) = w} h(i,\eta)$,
\begin{eqnarray*}
A^{-}_{i}(k) &= &\Phi \biggl( \frac{1 - \mu^{-}_{i}(k)}{\sigma
(i)} \biggr) -
\Phi \biggl( \frac{ - \mu^{-}_{i}(k)}{\sigma
(i)} \biggr),\\
 A^{+}_{i}(k)& = &\Phi
\biggl( \frac{1 - \mu^{+}_{i}(k)}{\sigma
(i)} \biggr) - \Phi \biggl( \frac{0 - \mu^{+}_{i}(k)}{\sqrt{\sigma(i)}} \biggr)
\end{eqnarray*}
and
\[
x_i\bigl(k,\sigma(i)\bigr) = \frac{\mu^{-} + \mu^{+}}{2} - \frac{\sigma
(i)^2}{\mu^{+}- \mu^{-}}
\log \frac{A^{-}_{i}(k)}{A^{+}_{i}(k)}.
\]

As a second example we show that the decomposition presented in Theorem
\ref{theodecomp} can be effectively implemented in Gibbsian systems
with compact-valued spins. We take $ A= [-1,1]$ and introduce the
following definitions.

%
\begin{defin}
A pairwise potential is a collection $\{ J(i,j) ,( i, j) \in\Z^d
\times\Z^d \} $ of real numbers which satisfies
%
\begin{equation}
\label{BB3} J(i,i)=0, \sup_{i \in\Z^d} \sum_{j \in\Z^d}
\bigl|J(i,j)\bigr| <\infty.
\end{equation}
\end{defin}

In what follows we use the notation
\[
\Sigma_i = \sum_{ j \in\Z^d } \bigl|J(i,j)\bigr| .
\]

For any $ i \in\Z^d, $ let $\eta(i)$ be the value of the spin at site
$i$ in the configuration $ \eta\in S.$

%
\begin{defin}
A probability measure $\mu$ on $(S, \mathcal{S})$ is said to be a
Gibbs state
relative to the potential $\{J (i,j)\}$ if for all\vadjust{\goodbreak} $i \in\Z^d$, a
version of
the conditional probability density of $\eta(i),$ given $\eta(j), j
\neq i,$ is given by
\[
\mu\bigl( \eta(i) = a | \eta(j) \mbox{ for all } j \neq i  \bigr) =
\frac{\exp ( a \sum_{ j \neq i } J(i,j)\eta(j)  )
}{Z^\eta},
\]
where
\[
Z^\eta= \int_A \exp \biggl( a \sum
_{ j \neq i } J(i,j) \eta(j) \biggr) \varrho(da).
\]
\end{defin}

In the following we consider the interaction $J_\beta= \beta J ,$
where $\beta$ is a positive parameter.
The associated Gibbs measure will be denoted $\mu$ without indicating
explicitly the dependence on $\beta.$
Now, put
%
\begin{equation}
\label{A2}c_i (a, \eta) = e^{\beta a \sum_{ j \in \Z^d }J(i,j)
\eta(j) } .
\end{equation}
Then, by construction, the process $(\sigma_t)$ with generator (\ref
{eqgenerator})
and this choice of change rates is reversible with
respect to the Gibbs state $\mu$ corresponding to the potential $
J_\beta
(i,j) = \beta J(i,j) .$ It is immediate to see that condition (\ref
{BB3}) implies the continuity condition
(\ref{eqcontinuity}).

We now give the explicit decomposition in one specific case.

\begin{ex}\label{exam2}
The following example is a Gibbsian time evolution with infinite range
interaction. The decomposition we present here is inspired by the one
presented in Galves, L\"ocherbach and Orlandi (\citeyear{GalLocOrl10}) in the case of
two color systems. In Galves, L\"ocherbach and Orlandi (\citeyear{GalLocOrl10}), for coupling reasons, it
was convenient to give a slightly different decomposition. The goal
there was to be able to couple together the infinite range Gibbsian
system with the finite range Gibbsian system obtained by truncating
the potential interaction. We suppose that the spin distribution
$\varrho$ is symmetric.
\end{ex}

Define for any $i \in\Z^d$ and any $ k
\geq- 1 , $
\[
S_i^{> k } := \sum_{j : \| i - j \| > k } \bigl|
J(i,j) \bigr| ,  \qquad S_i^{\le k } := \sum_{j : \| i - j \| \le k }
\bigl| J(i,j) \bigr| .
\]
Note that $ \Sigma_i = S_i^{ > - 1} .$

Then the decomposition (\ref{cmmc2}) holds with
%
\begin{equation}
\label{eqmi} M_i = \int_0^1 \bigl(
e^{ a \beta\Sigma_i } + e^{ - a \beta\Sigma
_i} \bigr) \varrho(da) .
\end{equation}
Moreover,
%
\begin{equation}
\label{eqalphafirst} \alpha_i ( -1) = 2 \int_0^1
e^{ - a \beta\Sigma_i } \varrho(da)
\end{equation}
and
%
\begin{equation}
\label{eqalphasecond} \alpha_i (k) = M_i + \int
_0^1 e^{ a \beta S_i^{ \le k } } e^{ - a
\beta S_i^{> k} }
\varrho(da) - \int_0^1 e^{a \beta\Sigma_i }
\varrho(da) .\vadjust{\goodbreak}
\end{equation}
Finally,
%
\begin{equation}
\label{eqlambdafirst} \lambda_i ( -1) = 2\frac{ \int_0^1 e^{ - a \beta\Sigma_i } \varrho
(da) }{\int_0^1  ( e^{ a \beta\Sigma_i } + e^{ - a \beta\Sigma
_i}  ) \varrho(da)}
\end{equation}
and
%
\begin{eqnarray}
\label{eqlambdasecond}
\quad \lambda_i(k)
&=&\int_0^1 e^{ a \beta S_i^{ \le k-1} } e^{ - a \beta S_i^{ > k
}} \nonumber\\[-2pt]
&&\hspace*{13pt}{}\times ( e^{ a \beta\sum_{ j : \| j-i \| = k } | J(i,j)|} -
e^{ -a \beta\sum_{ j : \| j-i \| = k } | J(i,j)|} ) \varrho
(da)
\\[-2pt]
&&{}\Big/ \int_0^1  \bigl( e^{ a \beta\Sigma_i } + e^{ - a \beta
\Sigma_i}  \bigr) \varrho(da) .\nonumber
\end{eqnarray}

\section{Perfect simulation}\label{perfect}
The goal of this section is to give an application of the Kalikow-type
decomposition given by Theorem~\ref{theodecomp}. This application is
a perfect simulation algorithm for the invariant measure of an
interacting multicolor system. We assume that the interaction rates
are continuous in the sense of (\ref{eqcontinuity}) and satisfy a
high noise condition. The basis of
the algorithm is the convex decomposition given in Theorem
\ref{theodecomp}. First of all, Proposition~\ref{theonstop}
gives a sufficient condition for exponential ergodicity which is based
on the construction of a dominating branching process.

From now on we will denote by $(\sigma^{\eta}_t)$ [and
$(\sigma^{\mu}_t)$] the multicolor system having generator $\LL$ given
by (\ref{eqgeneratorbis}) with a fixed initial configuration $\eta$
(a random configuration chosen with probability distribution
$\mu$).

%
\begin{prop}\label{theonstop}
Let $(c_i)_{i \in\Z^d} $ be a family
of rate functions satisfying the conditions of Theorem
$\ref{theodecomp}$. Furthermore, assume that
%
\begin{equation}
\label{eqm0} \underbar M = \inf_{ i \in\Z^d} M_i > 0
\end{equation}
and
%
\begin{equation}
\label{eqcondition2} \sup_{i \in\Z^d} \sum_{k \ge0}
\bigl|V_i (k)\bigr| \lambda_i (k) = \gamma < 1.
\end{equation}
Then the following two statements hold.
\begin{longlist}[(1)]
\item[(1)]
The process $(\sigma_t)$ admits a unique invariant probability measure
$\mu.$
\item[(2)]
For any finite set of sites $F \subset\Z^d ,$
for any $T > 0$ and any
initial configuration $ \eta,$ there exists a
coupling between the process $(\sigma^\eta_t)$ and the stationary
process $(\sigma^\mu_t)$
such that
\[
P\bigl( \sigma_T^\eta(F) \neq\sigma^\mu_T
(F) \bigr) \le|F| e^{- \underline
{ M} (1 - \gamma) T } .
\]
\end{longlist}
\end{prop}

Let us compare the above proposition to known results in the literature
on particle systems.
\begin{longlist}
\item[(1)] Condition (\ref{eqcondition2}) is stronger than Liggett's
existence condition (\ref{A3}) which does not imply the uniqueness of
the invariant measure.\vadjust{\goodbreak}
\item[(2)]
Let us compare our result to the $M < \varepsilon$-criterion of
Theorem 4.1 of
\citet{Lig85N2}, page 31. Recall that Liggett's quantity $M$ (translated
into our context) is given by
\[
M = \sup_{i \in\Z^d } \sum_{ j \neq i } \sup \bigl\{ \bigl\|
c_i ( a, \eta) \varrho(da ) - c_i ( a, \zeta) \varrho(da)
\bigr\|_{\mathrm{TV}} \dvtx \eta(k) = \zeta(k )\ \forall k \neq j \bigr\} .
\]
By our decomposition of $c_i ( a, \eta),$ this expression can be upper
bounded by
\[
M \le\sup_{i \in\Z^d} M_i \sum_{k \geq0 }
\lambda_i (k) \bigl|V_i (k) \bigr| .
\]
Since $\sup_i M_i < \infty, $ condition $ \sup_{i \in\Z^d} \sum_{k \ge0}   |V_i (k)| \lambda_i (k) < \infty$
implies condition~(3.8) of \citet{Lig85N2}.

Concerning the quantity $\varepsilon$ defined on page 24 of \citet{Lig85N2}, note that in our case, it can be written as follows:
\[
\varepsilon= \inf_{ i \in\Z^d} \inf_{ \eta,   a \neq b } \bigl[c_i
\bigl( a, \eta^{i,b}\bigr) \varrho\bigl( \{ a\}\bigr) + c_i
\bigl( b , \eta^{i,a} \bigr) \varrho\bigl( \{ b\}\bigr) \bigr],
\]
which will be equal to zero in general. Hence, with our techniques we
are able to treat cases where the $M< \varepsilon$-condition of
\citet{Lig85N2}, Theorem~4.1, is not satisfied.

\item[(3)] Condition (\ref{eqcondition2}) is a high-noise condition
reminiscent of Dobrushin--Shlosman condition [see \citet{MaeShl91}]. It is a sufficient condition ensuring that there is no phase
transition.
\end{longlist}

We are now in position to present the perfect simulation
scheme. Suppose we want to sample the configuration at site $ i $
under the invariant measure~$\mu.$ In a first step, we determine the
set of sites whose
spins influence the spin at site $i $ under equilibrium. We call this
set of sites \textit{ancestors} of $i $ and this stage \textit{backward
sketch procedure}. First, we climb up from time $0$ using a reverse
time Poisson point process with rate $M_i$. We stop when the last
Poisson clock before time $0$ rings. At that time, we choose a range
$k$ with probability $ \lambda_i (k)$. If $ k = - 1,$ we decide the
value of the spin using the law $ p_i^{ [ - 1]} \,d\varrho,$
independently of everything else. If $k $ is different from $-1,$ we
restart the above procedure from every site $ j \in V_i (k).$ The
procedure stops once each site involved has chosen range $ - 1.$ When
this occurs, we can start the second stage, in which we go back to the
future assigning spins to all sites visited during the first stage. We
call this procedure forward spin assignment procedure. This is done
from the past to the future by using the update probability densities
$p_i^ {[k]} $ starting at the sites which ended the first procedure by
choosing range $-1.$ For each one of these sites a spin is chosen
according to $p^{ [ - 1]} \,d\varrho.$ The values obtained in this way
enter successively in the choice of the values of the spins depending
on a neighborhood of range greater or equal to $0.$

We now give the precise form of the algorithm. Fix a finite set $F
\subset\Z^d .$ The following variables will be used:\vadjust{\goodbreak}

\begin{itemize}
\item$N$ is an auxiliary variable taking values in the set of
nonnegative integers $ \{ 0, 1,2, \ldots\} $.
\item$N^{(F)}_{\mathrm{STOP}}$ is a counter taking values in the set of
nonnegative integers $ \{ 0, 1, 2, \ldots\} $.
\item
$I $ is a variable taking values in $\Z^d$.
\item
$K$ is a variable taking values in $\{ -1, 0, 1, \ldots\}$.
\item
$B $ is an array of
elements of $\Z^d \times\{ -1, 0 , 1 , \ldots\} $.
\item
$C$ is a variable taking values in the set of finite subsets of $\Z^d$.
\item$ W $ is an auxiliary variable taking values in $ A^*$.
\item$\sigma$ is a function from $\Z^d $ to $A^*$.
\end{itemize}

\begin{algo}[(Backward sketch procedure)]\label{alg1}

\begin{enumerate}[10.]
\item[1.]
\textit{Input:} $F $; \textit{Output:} $N^{(F)}_{\mathrm{STOP}}$, $B$
\item[2.]
$N \leftarrow0,$ $N^{(F)}_{\mathrm{STOP}} \leftarrow0 ,$ $ B \leftarrow
\varnothing,$
$ C \leftarrow F $
\item[3.]
WHILE {$C \neq\varnothing$ }
\item[4.]
$N \leftarrow N+1 $
\item[5.] Choose randomly a position $I \in C$ and an integer $K\ge-1$
according to the probability distribution
\[
P(I=i, K=k )=\frac{M_i\lambda_i(k)}{ \sum_{j \in C}\sum_{l \ge
-1}M_j\lambda_j(l)}
\]
\item[6.] IF {$K = -1,$} {$C \leftarrow C \setminus\{ I \}$}
\item[7.]
ELSE $ C \leftarrow C \cup B_I (K)$
\item[8.]
ENDIF
\item[9.]$ B (N) \leftarrow (I, K ) $
\item[10.]
ENDWHILE
\item[11.]$N^{(F)}_{\mathrm{STOP}} \leftarrow N $
\item[12.]
RETURN $N^{(F)}_{\mathrm{STOP}}$, $B$
\end{enumerate}
\end{algo}

Now we use the following forward spin assignment
procedure to sample from the invariant measure $\mu .$
Recall that the choice of $\Delta$ in (\ref{eqgeneratorbis})
implies that the system does not change its colors. This explains
Step 9 in Algorithm~\ref{alg2}.

\begin{algo}[(Forward spin assignment procedure)]\label{alg2}
\begin{enumerate}[10.]
\item\textit{Input:} $N^{(F)}_{\mathrm{STOP}}$, $B$; \textit{Output:} $\{(i,\sigma
(i))\dvtx  i \in F \}$
\item$ N \leftarrow N^{(F)}_{\mathrm{STOP} } $
\item$\sigma(j) \leftarrow\Delta$ for all $ j \in \Z^d $
\item WHILE {$N \ge1$}
\item$ (I,K) \leftarrow B(N)$
\item IF {$K= -1 $} choose $W $ randomly in $A$
according to the probability distribution
\[
p_I^{[-1]} \,d\varrho
\]
\item ELSE {choose $W $ randomly in $A^*$
according to the probability distribution
\[
p_I^{[K]} ( \cdot| \sigma) \,d\varrho^*
\]
}
\item ENDIF
\item IF {$ W \neq\Delta $} put $ \sigma(I) \leftarrow W $
\item ENDIF
\item{ $ N \leftarrow N-1 $}
\item ENDWHILE
\item RETURN $ \{ (i, \sigma(i)) \dvtx  i \in F \} $
\end{enumerate}
\end{algo}

The next theorems summarize the properties of Algorithms~\ref{alg1} and~\ref{alg2}.
In order to distinguish clearly to which part of the two algorithms we
refer, we shall write
$P_{\mathrm{sketch} } $ for the probability associated to the backward sketch procedure.

%
\begin{theo}\label{theo15}
Suppose that
%
\begin{equation}
\label{eqexactcondition} \mbox{for all $ F \subset\Z^d $ finite}\qquad
P_{\mathrm{sketch}} \bigl( N^{(F)}_{\mathrm{STOP}} < \infty\bigr) = 1 .
\end{equation}
Then the following two statements hold.
\begin{longlist}[(1)]
\item[(1)]
Algorithms~\ref{alg1} and~\ref{alg2} are successful.
\item[(2)]
The process $ (\sigma_t)$ admits a unique invariant measure $\mu.$
The law of the set
$\{ (i, \sigma(i)) \dvtx  i \in F\} $ printed at the end of Algorithms~\ref{alg1}
and~\ref{alg2} is the projection on
$A^F$ of $\mu.$
\end{longlist}
\end{theo}

The next theorem states sufficient conditions ensuring (\ref
{eqexactcondition})
and gives also a control on the rate of convergence.

%
\begin{theo}\label{theo2}
\textup{(1)} The sub-criticality condition $(\ref{eqcondition2})$ implies $(\ref
{eqexactcondition})$. More
precisely, we have
%
\begin{equation}
\label{s4} P_{\mathrm{sketch}} \bigl( N^{(F)}_{\mathrm{STOP}} > N \bigr)
\le|F| \gamma^N ,
\end{equation}
where $\gamma$ is given in $(\ref{eqcondition2})$.\vspace*{-6pt}
\begin{longlist}[(2)]
\item[(2)]
Suppose, in addition to $(\ref{eqcondition2})$, that $(\ref{eqm0})$
holds. Fix a time $t > 0,$ some finite set of sites $F \subset\Z^d $
and two
initial configurations $ \eta$ and $\zeta\in A^{\Z^d} .$ Then there
exists a
coupling of the two processes $(\sigma^\eta_s)_s$ and $(\sigma^\zeta_s)_s $ such that
\[
P\bigl( \sigma_t^\eta(F) \neq\sigma^\zeta_t
(F) \bigr) \le|F| e^{-
\underline{M} ( 1 - \gamma) t } .
\]
\end{longlist}
\end{theo}

The proofs of Proposition~\ref{theonstop}, Theorems~\ref{theo15} and~\ref{theo2}
will be given in the next section.

\section{\texorpdfstring{Proofs of Proposition \protect\ref{theonstop},
Theorems \protect\ref{theo15} and \protect\ref{theo2}}
{Proofs of Proposition 1, Theorems 2 and 3}} \label{secbw}

The proofs rely on the notion of
\textit{black and white time-reverse sketch process} that we will
introduce now. The black and white time-reverse sketch process gives the\vadjust{\goodbreak}
mathematically precise description of the backward black and white
Algorithm~\ref{alg1} given in Section~\ref{perfect}.

For each $ i \in\Z^d, $ denote
by $\cdots T_{-2}^i <T_{-1}^i < T_{0}^i < 0 < T_1^i < T_2^i <
\cdots$ the occurrence times of the rate $M_i$ Poisson point process
$N^i $ on the real line. The Poisson point processes associated to
different sites are independent. To each point $T_n^i$ associate an
independent mark $K^i_n$ according to the probability distribution
$(\lambda_i(k))_{k \ge-1}$. As usual, we identify the Poisson point
processes and the associated counting measures.

For each $i \in\Z^d$ and $t \in\R$ we define the time-reverse point
process starting at time $t,$ associated to site $i,$
%
\begin{eqnarray}
\label{eqtildet} \tilde{T}^{(i,t)}_n &=& t -
T^i_{N^i(0,t]-n+1},\qquad t \ge 0,
\nonumber
\\[-8pt]
\\[-8pt]
\nonumber
\tilde{T}^{(i,t)}_n & = & t - T^i_{-N^i(t,0]-n+1},\qquad
t < 0 .
\end{eqnarray}
To these time-reverse point processes we can associate in an obvious
way the corresponding marks $ \tilde{K}^{(i,t)}_n, n \in\Z .$
Finally, for each site $i \in\Z^d$, $k \ge-1$, the reversed $k$-marked
Poisson point process returning from time $t$ is defined as
%
\begin{equation}
\label{eqtilden} \tilde{N}^{(i,t,k)}[s,u] = \sum
_{n} \one_{\{s \le\tilde
{T}^{(i,t)}_{n} \le u\}} \one_{\{\tilde{K}^{(i,t)}_n = k\}}.
\end{equation}

To define the black and white time-reverse sketch process we need to
introduce a family of transformations $\{\pi^{(i,k)}, i \in\Z^d, k
\ge-1 \}$ on the set of finite subsets of $\Z^d,$ $ \mathcal{F}(\Z^d),$
defined as follows. For any unitary set $\{j\}$,
%
\begin{equation}\label{eqpij}
\pi^{(i,k)}\bigl(\{j\}\bigr)= \left\{
\begin{array}{l}
V_i (k),\qquad\mbox{if } j=i
\\
\{j\},\qquad\hspace*{10pt} \mbox{otherwise}
\end{array}
\right\} .
\end{equation}
Notice that for $k=-1,$ $ \pi^{(i,k)}(\{i\}) = \varnothing.$ For any
finite set $F \subset\Z^d$, we similarly define
%
\begin{equation}
\label{eqpif} \pi^{(i,k)}(F) = \bigcup_{j \in F}
\pi^{(i,k)}\bigl(\{j\}\bigr) .
\end{equation}

The black and white time-reverse sketch process starting at site $i$
at time $t$ will be denoted by $(C_s^{(i,t)})_{s \geq0}.$
$C_s^{(i,t)}$ is the set of sites at time $s$ whose colors affect the
color of site $i$ at time $t.$ We call this set $C_s^{(i,t)} $ \textit{set
of ancestors
of $i$ at time $s$ before time $t.$} The evolution of this process is
defined through the following equation: $C_0^{(i,t)} := \{i\},$ and
%
\begin{eqnarray}
\label{eqct}\qquad f\bigl( C_s^{(i,t)}\bigr) &=& f
\bigl(C_0^{(i,t)}\bigr)
\nonumber
\\[-8pt]
\\[-8pt]
\nonumber
&&{}+ \sum_{k \ge-1}
\sum_{j
\in\Z^d} \int_0^s
\bigl[f\bigl(\pi^{(j,k)} \bigl(C_{u-}^{(i,t)}\bigr)
\bigr) - f\bigl(C_{u-}^{(i,t)}\bigr)\bigr] \tilde{N}^{(j,t,k)}(du),
\end{eqnarray}
where $f\dvtx  \mathcal{F}(\Z^d) \rightarrow\R$ is any bounded cylindrical
function. This family of equations characterizes completely the time\vadjust{\goodbreak}
evolution $\{C_s^{(i,t)}, s \ge0\}$. For any finite set $F \subset
\Z^d$ define
\[
C_s^{(F,t)} = \bigcup_{i \in F}
C_s^{(i,t)}.
\]

The following proposition summarizes the properties of the family of
processes defined above.

%
\begin{prop}
For any finite set $F \subset\Z^d$, $\{C_s^{(F,t)}, s \ge0\} $ is a
Markov jump
process having as infinitesimal generator
%
\begin{eqnarray}
\label{eqgeneratord} L f(C) &=& M_i \sum_{i \in C}
\sum_{k \ge0} \lambda_i (k) \bigl[f\bigl(C
\cup V_i(k)\bigr) - f(C)\bigr]
\nonumber
\\[-8pt]
\\[-8pt]
\nonumber
 &&\hspace*{46pt}{}+ \lambda_i (-1) \bigl[f
\bigl(C \setminus\{i\}\bigr) - f(C)\bigr] ,
\end{eqnarray}
where $f$ is any bounded cylindrical function.
\end{prop}

\begin{pf} The proof follows in a standard way from the construction
\reff{eqct}.
\end{pf}

If we are interested in simulating from the invariant measure of the
process, then we will start the black and white time-reverse sketch
process at time $ t = 0 ;$ if, however, we wish to construct the process
at time $t,$ we shall start the black and white time-reverse sketch
process at that time $t $ precisely.

\subsection{Backward oriented percolation and sub-criticality}
For the algorithm to be successful it is crucial to show that
$\bigcup_{ s \geq0 } C_s^{(i,t)}$, the set of ancestors of any site~$i$,
is finite with probability one. Formally, let
\[
T^{(i)}_{\mathrm{STOP}} = \inf\bigl\{ s \dvtx C_s^{(i,0)}
= \varnothing\bigr\}
\]
be the relaxation time.
We
introduce the sequence of successive jump times $\tilde T_n^{(i)} , n
\geq1 ,$ of processes $N^{(j,k)} $ whose jumps occur in
(\ref{eqct}), for $t = 0.$ Let $\tilde T_1^{(i)} = T_1^{(i,0)} $ and define
successively for $n\geq2 $
%
\begin{equation}
\label{Ajumps} \tilde T_n^{(i)} = \inf\bigl\{ t > \tilde
T_{n-1}^{(i)} \dvtx \exists j \in C^{(i, 0)}_{\tilde T_{n-1}^{(i)} }
, \exists k \dvtx N^{(j,k)} \bigl( \bigl]\tilde T_{n-1}^{(i)}
, t \bigr]\bigr) = 1 \bigr\} .
\end{equation}
We write $\tilde K_n^{(i)} $ for the associated marks.
Now we put
%
\begin{equation}
\label{RR5} \mathbf{C}^{(i)}_n = C^{(i,0)}_{ \tilde T_n^{(i)} }
\end{equation}
and
\[
N^{(i)}_{\mathrm{STOP}} = \inf\bigl\{ n \dvtx \mathbf{C}^{(i)}_n
= \varnothing\bigr\} .
\]
This is the number of steps of the backward sketch process---and it
is exactly the number of steps of Algorithm~\ref{alg1}. For the perfect
simulation algorithm to be successful, it is crucial to show that 
the number of steps $
N^{(i)}_{\mathrm{STOP}} $ is finite. Since at every step of the algorithm a
finite interaction range\vadjust{\goodbreak}
$k$ is chosen, this implies automatically that also $ T^{(i)}_{\mathrm{STOP}}$
is finite almost surely.
However, in order to control the speed of convergence, we need a
precise control on the tail
probabilities of $ T^{(i)}_{\mathrm{STOP}}.$ To this aim we estimate the
volume of the set $ C_s^{(F, t)}= \bigcup_{i \in F} C_s^{(i,t)}$ where~$F$ is a bounded set of $\Z^d$.
%
\begin{lem}\label{lemma1E}
%
\begin{equation}
\label{eqgron2} E \bigl(\bigl |C_s^{(F, t)} \bigr|\bigr) \le|F|
e^{ - \underline{M} ( 1 - \gamma) s,}
\end{equation}
where $\underline{M}$ is defined in \eqref{eqm0} and $\gamma$ in
\eqref{eqcondition2}.
\end{lem}
\begin{pf}
Fix some $N \in\N.$ Let $L^i_s = | C_s^{(i,t)}| $ and
\[
T_N = \inf\bigl\{ t \dvtx L^i_t \geq N
\bigr\} .
\]
Then by (\ref{eqct}),
%
\begin{eqnarray}
\label{equb} L^i_{s \wedge T_N} & \le& 1 + \sum
_{k \geq1} \sum_{j \in\Z^d} \int
_0^{s \wedge T_N} \bigl[\bigl|V_j(k)\bigr| - 1\bigr]
1_{\{ j \in C^{(i,t)}_{u-} \}} \tilde{N}^{(j,t,k)}(du)
\nonumber
\\[-8pt]
\\[-8pt]
\nonumber
&&{} -\sum_{j \in\Z^d} \int_0^{s \wedge T_N}
1_{\{ j \in
C^{(i,t)}_{u-} \}} \tilde{N}^{(j,t,-1 )}(du) .
\end{eqnarray}
Recall that $\underline{M} = \inf_{i \in\Z^d} M_i > 0.$
Passing to expectation and using that, by condition (\ref{eqcondition2}),
\[
M_j \biggl( \biggl(\sum_{k \geq1}
\lambda_j (k) \bigl[ \bigl|V_j(k)\bigr| - 1\bigr] \biggr) -
\lambda_j (-1) \biggr) \le- \underline{M} ( 1 - \gamma) < 0 ,
\]
which yields
%
\begin{eqnarray}
\label{equpperbound2} E \bigl(L^i_{s \wedge T_N}\bigr)& \le& 1 +
\sum_{j \in\Z^d} M_j \biggl( \biggl(\sum
_{k \geq1} \lambda_j (k) \bigl[
\bigl|V_j(k)\bigr| - 1\bigr] \biggr) - \lambda_j (-1) \biggr)
\nonumber
\\
&&\hspace*{38pt}{} \times E \int_0^{s\wedge T_N} 1_{\{ j \in C^{(i,t)}_{u-} \}} \,du
\\
&\le& 1 - \underline{M} ( 1 - \gamma) E \int_0^{s\wedge T_N}
L^i_u\, du .\nonumber
\end{eqnarray}
Letting $N \to\infty,$ we thus get by Fatou's lemma that
\[
E\bigl(L^i_s\bigr) \le1 - \underline{M} (1 - \gamma)
\int_0^s E\bigl(L^i_u
\bigr) \,du .
\]
This implies that
\[
E\bigl(L^i_s\bigr) \le1 \qquad\mbox{for all } s \geq0 .
\]
Hence, we may apply Gronwall's lemma which yields
%
\begin{equation}
\label{eqgron} E\bigl(L^i_s\bigr) \le e^{ - \underline{M} (1 - \gamma) s} .\vadjust{\goodbreak}
\end{equation}
Hence, since $ |C_s^{(F, t)} | \le\sum_{i \in F } |C_s^{(i, t)}| =
\sum_{i \in F } L_s^i,$
%
\begin{equation}
\label{eqgron2} E \bigl(\bigl |C_s^{(F, t)} \bigr|\bigr) \le|F|
e^{ - \underline{M} ( 1 - \gamma) s .}
\end{equation}
\upqed\end{pf}

\subsection{\texorpdfstring{Proof of Proposition \protect\ref{theonstop} and
Theorem \protect\ref{theo15}}
{Proof of Proposition 1 and Theorem 2}}
Proposition~\ref{theonstop} is an immediate consequence of Theorem
\ref{theo15}, item 2, and Theorem~\ref{theo2}, item 2.

\begin{pf*}{Proof of Theorem~\ref{theo15}}
Item 1 of Theorem~\ref{theo15} is evident. We give the proof of item 2
of Theorem~\ref{theo15}.

Write $\mu_F$ for the law of the output $\{ (i,\sigma(i)) \dvtx  i \in F \}
$ of Algorithms~\ref{alg1} and~\ref{alg2}; $\mu_F$ is a probability measure on $ ( A^F,
\mathcal{A}^F ).$ By construction, the family of probability laws $ \{
\mu_F,   F \subset\Z^d \mbox{ finite} \} $ is a consistent family
of finite dimensional distributions. Hence, there exists a unique
probability measure $ \mu$ on $ ( S, \mathcal{S})$ such that $ \mu_F
$ is the projection onto $A^F$ of $\mu, $ for any fixed finite set of
sites $F \subset\Z^d.$

We show that $\mu$ is the unique invariant measure of the process
$(\sigma_t) .$
In order to do so, we use a slight modification of Algorithms~\ref{alg1} and~\ref{alg2}
in order to
construct $\sigma_t^\eta,$ for some fixed initial configuration $\eta
\in S.$ The modification is defined
as follows. Let $ T_{\mathrm{STOP}} $ and $T$ be variables taking values in $(0,
\infty) .$
Replace Steps 1--3 of Algorithm~\ref{alg1} by:

\noindent 1. Input: $ F ; \mbox{ Output: } N^{(F)}_{\mathrm{STOP}}, B, C$

\noindent 2. $N \leftarrow0,$ $N^F_{\mathrm{STOP}} \leftarrow0 ,$ $ B \leftarrow
\varnothing,$ $ C \leftarrow F, $ $T_{\mathrm{STOP}} \leftarrow0 $

\noindent 3. WHILE $ T_{\mathrm{STOP}} < t $ and $ C \neq\varnothing$

\noindent 3$^{\prime}$. Choose a time $T \in(0, +\infty) $ randomly according to the
exponential distribution with parameter
$ \sum_{j \in C} M_j .$ Update
\[
T_{\mathrm{STOP}} \leftarrow T_{\mathrm{STOP}} + T .
\]

Finally replace Step 12 of Algorithm~\ref{alg1} by:

\noindent 12. RETURN $ N^{(F)}_{\mathrm{STOP}}, B, C$.

In this modified version, we stop the algorithm after time $t,$ hence,
the output
set $C$ might not be empty. The output $C$ is exactly the set $C_t^{(F,
t)}, $ the set of
sites at time $0$ whose colors influence the colors of sites in $F$ at
time $t.$ Finally, notice
that if $C = \varnothing, $ then $ T_{\mathrm{STOP}} < t ,$ and in this case, $
T_{\mathrm{STOP}} = T^{F}_{\mathrm{STOP}} $ is the relaxation time introduced in the
previous subsection.

Concerning Algorithm~\ref{alg2}, replace Step 1 of Algorithm~\ref{alg2} by:

\noindent 1. Input: $ N^{(F)}_{\mathrm{STOP}}, B, C; $ Output: $ \{ (i, \sigma(i) ) : i
\in F\}$

and Step 3 by:

\noindent 3. $\sigma(j) \leftarrow\eta(j) $ for all $j \in C;$ $\sigma(j)
\leftarrow\Delta$ for all $ j \in \Z^d \setminus C .$

Then the law of the set $\{ (i, \sigma(i) ) \dvtx  i \in F \} $ printed at
the end of the modified Algorithm~\ref{alg2} is the law of $ \sigma^\eta_t (F)
.$ Notice that the output of the modified Algorithm~\ref{alg2} equals the output
of the unmodified Algorithm~\ref{alg2} if $ T_{\mathrm{STOP}} < t . $

We first give an intuitive argument showing that $\mu$ must be
invariant for $(\sigma_t).$ Write $P_t$ for the
transition semigroup of $(\sigma_t).$ Fix
$t> 0 $ and a finite set of sites $F \subset\Z^d.$ Suppose we want to
determine the projection on $A^F$ of $\mu P_t .$
This means that
we have first to run the above modified Algorithm~\ref{alg1} up to time $t.$ It
gives as output the set of ancestor sites $C= C_t^{(F,t)}.$
We then have to run the modified Algorithm~\ref{alg2} with initial configuration
$ \{\sigma(i ), i \in C = C_t^{(F,t) } \}$ chosen according to $\mu$;
cf. Step 3. But this means that we have to concatenate the modified
Algorithm~\ref{alg1} with the original Algorithm~\ref{alg1}, where the Algorithm~\ref{alg1} is now
starting from $ C_t^{(F,t)} $ instead of $F.$ In other words, we
concatenate two backward sketch processes and consider $C^{(C^{(F,t)},
0)} $ up to the time of its extinction. By the Markov property and the
stationarity of the Poisson processes, this means that we directly
consider $C^{(F,t)} $ up to the time of its extinction---which is
finite by our assumptions. Hence, $ \mu P_t = \mu$ in restriction to
finite cylinder sets.

We now give a more formal argument. Recall that the output of the
modified Algorithm~\ref{alg2} equals the output
of the unmodified Algorithm~\ref{alg2} if $ T_{\mathrm{STOP}} < t . $ Let $ f\dvtx  A^F \to\R_+ $ be a bounded measurable function.
Then
%
\begin{eqnarray}
E \bigl[ f \bigl( \sigma_t^\eta(i), i \in F \bigr) \bigr]
&= &E \bigl[ f \bigl( \sigma_t^\eta(i), i \in F \bigr),
T_{\mathrm{STOP}} < t \bigr] \nonumber\\
&&{}+ E \bigl[ f \bigl( \sigma_t^\eta(i),
i \in F \bigr) , T_{\mathrm{STOP}} \geq t \bigr]
\nonumber
\\[-8pt]
\\[-8pt]
\nonumber
& = & E \bigl[ f \bigl( \sigma(i), i \in F \bigr), T^F_{\mathrm{STOP}}
< t \bigr]
\\
&&{} + E \bigl[ f \bigl( \sigma_t^\eta(i), i \in F \bigr) ,
T_{\mathrm{STOP}} \geq t \bigr] ,\nonumber
\end{eqnarray}
where $ ( \sigma(i), i \in F)$ is the output of the unmodified
Algorithms~\ref{alg1} and~\ref{alg2}.

But
\[
E \bigl[ f \bigl( \sigma_t^\eta(i), i \in F \bigr) ,
T_{\mathrm{STOP}} \geq t \bigr] \le\| f \|_{\infty} P_{\mathrm{sketch}}
\bigl( T^F_{\mathrm{STOP}} \geq t\bigr) \to0\qquad \mbox{as } t \to \infty,
\]
since finiteness of $N_{\mathrm{STOP}}^F $ implies the finiteness of $T^F_{\mathrm{STOP}}
.$ Hence, we obtain that
\[
\lim_{t \to\infty} E \bigl[ f \bigl( \sigma_t^\eta(i),
i \in F \bigr) \bigr] = E \bigl[ f \bigl( \sigma(i), i \in F \bigr) \bigr] ,
\]
since $ 1_{\{ T^F_{\mathrm{STOP}} < t \} } \to1 $ almost surely.

This implies that $\mu$ is an invariant measure of the process.
Replacing the initial condition $\eta$ by any stationary initial
condition, we finally also get uniqueness of the invariant measure.
Thus Theorem~\ref{theo15} is proved.
\end{pf*}

\subsection{\texorpdfstring{Proof of Theorem \protect\ref{theo2}}{Proof of Theorem 3}}

We start by proving (\ref{s4}). Let
\[
L^{(i)}_n = \bigl| \mathbf{C}^{(i)}_n\bigr|
\]
be the cardinal of the set $\mathbf{C}^{(i)}_n$ after $n$ steps of the
algorithm [recall (\ref{RR5})]. Then due to our assumptions, $
L^{(i)}_n$ can be compared to a multi-type branching process $Z_n$
having offspring mean which is bounded by $\gamma$ at each step, such
that $ L^{(i)}_n \le Z_n $ for all $n.$ The details are given in
Galves,\vadjust{\goodbreak} L\"ocherbach and Orlandi (\citeyear{GalLocOrl10}), proof of Theorem 1, and are
omitted here. Thus,
\[
P\bigl( N^{(i)}_{\mathrm{STOP}} > n \bigr) = P\bigl(
L^{(i)}_n > 0 \bigr) = P\bigl( L^{(i)}_n
\geq1 \bigr) \le P( Z_n \geq1 ) \le E (Z_n ) =
\gamma^n.
\]
When starting with the initial set $F$ instead of the singleton $\{i\},
$ then the above estimates remain true
by multiplying with $|F|,$ due to the independence properties of the
branching process.

Concerning item 2 of Theorem~\ref{theo2}, we use once more the
modified Algorithms~\ref{alg1} and~\ref{alg2} introduced in the proof of Theorem \ref
{theo15} above. In order to realize the coupling, we use the same realizations
of $T, I $ and $ K $ for the construction of $\sigma_t^\eta$ and
$\sigma_t^\zeta$. Write $L_s$ for the cardinal of $C_s^{(F,t)}.$
Clearly, both realizations of $\sigma_t^\eta$ and $\sigma_t^\zeta$
do not depend on the initial configuration $\eta,$ $\zeta$,
respectively, if the output $C$ of the modified Algorithm~\ref{alg1} is
void. Thus, by Lemma~\ref{lemma1E},
\begin{eqnarray*}
P \bigl( \sigma_t^\eta ( F) \neq\sigma_t^\zeta(F)
\bigr) & \le& P( T_{\mathrm{STOP}} \geq t )
\\
& =& P ( L_t \geq1 )
\\
& \le& E (L_t) \le|F| e^{- \underline{M} ( 1 - \gamma) t } .
\end{eqnarray*}
This implies that the convergence toward the unique
invariant measure takes place exponentially fast.
The proof of Theorem~\ref{theo2} is complete.

\section{Applications for perfect simulation}\label{sectionappli}

\subsection{Maximum likelihood estimation in Gibbs distributions}
Parameter estimation for Gibbs distributions in the infinite lattice
is usually based on the maximum likelihood approach [see, e.g.,
Gidas (\citeyear{Gid88}, \citeyear{Gid91})]. The maximum likelihood estimation is theoretically
well understood in this framework. \citet{Com92} proved the
consistency of the MLE for exponential families of Markov random
fields on the lattice. Also, in the case of no phase transition,
Jan{\v{z}}ura (\citeyear{Jan97}) proved asymptotic normality and efficiency of the MLE
inside the uniqueness region of the Gibbs distributions
considered. \citet{ComGid92} considered maximum likelihood
estimators for Markov random fields over $\mathbf{Z}^d$ from incomplete
data. They prove the strong consistency of maximum
likelihood estimators in this case. Their
results hold irrespective of the presence of long-range correlations
or nonanalytic behavior of the underlying quantities. The parameter
space is thereby allowed to be noncompact.

However, the numerical feasibility of the ML method is strongly
limited, due to the computation of the normalizing constant for each
relevant parameter, in particular, for each temperature. \citet{GeyTho92}
devised a rather ingenious method for this computation
based on an MCMC computation of the equilibrium distribution for a
fixed value of the parameter.\vadjust{\goodbreak}

Consider the family of probability densities with respect to a
reference measure~$\mu$ given by
%
\begin{equation}
\label{eqgibbs} f(x,\theta) = \frac{1}{Z(\theta)} \exp\bigl\langle T(x),\theta
\bigr\rangle,
\end{equation}
where $\langle T(x),\theta\rangle$ denotes the inner product between the canonical
parameter $\theta$ and the sufficient statistics $T(x)$ and
%
\begin{equation}
\label{eqpartition} Z(\theta) = \int\exp\bigl\langle T(x),\theta\bigr\rangle \,d
\mu(x).
\end{equation}

Denote by $P_{\psi}$ the measure having density $f(\cdot,\psi)$ with
respect to $\mu$. Then,
%
\begin{equation}
\label{eqpartition2} Z(\theta) = Z(\psi) \int\exp\bigl\langle T(x),\theta- \psi
\bigr\rangle \,dP_{\psi}(x).
\end{equation}

Therefore, if we have $X_1, X_2, \ldots$ i.i.d. random objects with
distribution $P_{\psi}$, we have that
%
\begin{equation}\quad
\label{eqdn} d_n(\theta) = \frac{1}{n} \sum
_{i=1}^n \exp\bigl\langle T(X_i),
\theta- \psi\bigr\rangle \rightarrow d(\theta) = \frac{Z(\theta)}{Z(\psi)} \qquad\mbox{almost
surely.}
\end{equation}

The maximum likelihood of $\theta$ can be taken as
%
\begin{equation}
\label{eqmle} \hat{\theta} = \operatorname{arg max} \log f(x, \theta) + \log Z(
\psi) = \operatorname{arg max} \bigl\langle T(x),\theta\bigr\rangle - \log d(
\theta),
\end{equation}
and its Monte Carlo approximant
%
\begin{equation}
\label{eqmleaprox} \hat{\theta}_n = \operatorname{arg max} \bigl
\langle T(x),\theta\bigr\rangle - \log d_n(\theta).
\end{equation}

Notice that if we can perfectly simulate from $P_{\psi}$, we have
trivially that $\hat{\theta}_n \rightarrow\hat{\theta}$ as $n
\rightarrow\infty$ along with the rate of convergence of such
convergence.

\subsection{\texorpdfstring{Attaining Ornstein's $\bar d$-distance for ordered pairs
of Ising probability distributions}
{Attaining Ornstein's d-distance for ordered pairs
of Ising probability distributions}}
In this section we show how to use the decomposition of Theorem
\ref{theodecomp} and the above perfect simulation algorithm in order
to construct an explicit
coupling attaining Ornstein's $\bar d$-distance for
two ordered Ising probability measures. Let $A:=\{-1,1\}$ and $S =
A^{\Z^d }.$

We consider a ferromagnetic pairwise interaction $J,$ that is, a
collection $
\{J(i,j),i\neq j,i,j\in{\mathbb{Z}}^{d}\}$ of positive real numbers
satisfying $
J(i,j)=J(j,i)$ for all $i,j\in{\mathbb{Z}}^{d}$ and for all $i \in\Z^d,$
%
\begin{equation}
\label{B3} J(i,i) = 0 , \qquad\sup_{i\in{\mathbb{Z}}^{d}}\sum_{j}J(i,j)<
\infty.
\end{equation}

Let $\{\tilde{J}(i,j),i\neq j,i,j\in{\mathbb{Z}}^{d}\}$ be another pairwise
interaction satisfying an analogous summability assumption such that
\[
J(i,j)\leq\tilde{J}(i,j)\qquad\mbox{for all }i,j\in{\mathbb{Z}}^{d} .
\]

Moreover, let $ \{ h_i , i \in\Z^d \}$ and $ \{ \tilde h_i , i \in\Z^d \} $
be two collections of positive real numbers representing an external
field such that
\[
h_i \leq\tilde h_i \qquad\mbox{for all } i \in
\Z^d , \sup_i \tilde h_i < \infty.
\]

Finally we suppose that,
%
\begin{equation}
\label{eqsufficient} \mbox{for all } i \in\Z^d ,\qquad  \sum
_{j } \bigl[ \tilde J (i,j) - J(i,j) \bigr] \le\tilde
h_i - h_i .
\end{equation}

Recall that a probability measure $\mu$ on $S$ is said to be a Gibbs
measure relative
to the interaction $J$ and the external field $ \{ h_i\} $ if for all
$i\in{\mathbb{Z}}^{d}$ and for any fixed $
\zeta\in S,$ a version of the conditional probability $\mu(\{\sigma
\dvtx \sigma(i)=\zeta(i)|\sigma(j)=\zeta(j)\mbox{ for all }\ j\neq i\}
)$ is
given by
%
\begin{eqnarray}\label{B4}
&&\mu\bigl(\bigl\{\sigma\dvtx \sigma(i)=\zeta(i)|\sigma(j)=\zeta(j)\mbox { for all
}\ j\neq i\bigr\}\bigr)
\nonumber
\\[-8pt]
\\[-8pt]
\nonumber
&&\qquad=\frac{1}{1+\exp(-2 \beta[\sum_{j}J(i,j)\zeta(i)\zeta(j) + h_i \zeta(i) ])
}.
\end{eqnarray}

The Gibbs measure $\tilde\mu$ associated to the interaction $\tilde
J$ and the external field $ \{ \tilde h_i \}$
is introduced analogously.

We consider a Glauber dynamics $(\sigma_t (i), i \in\Z^d , t
\in{\mathbb{R}} )$ taking values in $S= A^{\Z^d} $ and having $\mu
$ as reversible
measure. The process is defined by the rates
$c_i (\sigma),$ $i \in\Z^d, $ where $c_i (\sigma)$ is the rate at
which the spin $i$ flips (i.e., changes its sign) when
the system is in the configuration $\sigma .$ We take
%
\begin{equation}
\label{A2} c_i (\sigma) = \exp \biggl( -\beta\biggl[ \sum
_{ j } J (i,j) \sigma(i ) \sigma(j) + h_i
\sigma(i) \biggr] \biggr).
\end{equation}

By construction, the process $(\sigma_t)_t$ is reversible with respect to
the Gibbs measure $\mu.$ In the same way, we can define a Glauber
dynamics $(\tilde\sigma_t)_t $
reversible with respect to the Gibbs measure $\tilde\mu,$ associated
to the interaction $\tilde J$ and the external field $ \{ \tilde h_i\} .$

The main idea of our approach is a coupled construction of the
processes $(\sigma_t )$ and $(\tilde\sigma_t)$
which is order preserving. More precisely, let us write
\[
\sigma \le\tilde\sigma \quad\mbox{if and only if}\quad \sigma(i ) \le \tilde\sigma (i )\qquad
\mbox{ for all } i \in\Z^d .
\]

Write
\[
\mathcal{S}= \bigl\{ ( \sigma, \tilde\sigma ) \in A^{ \Z^d } \times
A^{\Z^d } \dvtx \sigma \le\tilde\sigma \bigr\} = \bigl\{ (-1, -1 ), (-1,
+1), (+1, +1 ) \bigr\}^{\Z^d} .
\]




We now describe the coupled time evolution of $\sigma_t $ and $\tilde
\sigma_t.$ Start with an ordered couple of initial configurations
$\eta\le
\tilde\eta$ at time $0.$ Let
\[
M_i= 2e^{\beta[ \sum_{j \in{\mathbb{Z}}^d } \tilde J (i,j) + \tilde
h_i ] } .
\]
For each site $i \in\Z^d ,$ consider a rate $M_i$
Poisson process $N^i .$ The Poisson processes corresponding to distinct
sites are independent. If at time $t,$ the Poisson\vadjust{\goodbreak} clock ate site $i$ rings,
then both processes try simultaneously to update their spin at site $i.$
Process $\sigma$ replaces spin $\sigma(i) $ by $- \sigma(i) $ with
probability
\[
\frac{c_i ( \sigma_t)}{M_i},
\]
and the process $\tilde\sigma$ replaces spin $\tilde\sigma(i) $ by $-
\tilde\sigma(i) $ with probability
\[
\frac{\tilde c_i ( \tilde\sigma_t)}{M_i} .
\]

Hence, we can introduce the following probability measures. For any
configuration $\sigma$ with $\sigma (i ) = - 1 $, we put
\[
p_i ( +1 | \sigma) = \frac{c_i ( \sigma)}{M_i} ,\qquad p_i ( - 1 |
\sigma ) = 1 - p_i ( +1 | \sigma) .
\]

In the same way, for any configuration $\sigma$ with $\sigma(i ) =
+1,$ we put
\[
p_i ( -1 | \sigma) = \frac{c_i ( \sigma)}{M_i} ,\qquad p_i ( + 1 |
\sigma ) = 1 - p_i ( -1 | \sigma) .
\]

The same definitions hold for $\tilde c_i $ with obvious modifications.
Then we have by construction
and thanks to condition (\ref{eqsufficient}) that
\[
p_i ( +1 | \sigma) \le\tilde p_i ( +1 | \tilde\sigma) \qquad
\mbox{whenever } \sigma\le\tilde\sigma.
\]

This stochastic order makes it possible to construct a coupled Glauber
dynamics $(\sigma_t, \tilde\sigma_t)_t$ taking values in the space
of ordered
configurations $
\mathcal{S } .$

At each jump time $t$ of one of the Poisson processes $N^i, $ the ordered
configuration $(\sigma_t, \tilde\sigma_t) $ is replaced at site $i$
by the
ordered pair
%
\begin{eqnarray}
\qquad&& (+1, +1 ) \qquad\mbox{with probability }
P_i \bigl( ( +1, +1 ) | ( \sigma_t, \tilde
\sigma_t) \bigr) = p_i ( + 1 | \sigma_t ) ,
\nonumber\\
&&( - 1, -1) \qquad
\mbox{with probability } P_i \bigl( ( -1, -1 ) | (
\sigma_t, \tilde \sigma_t) \bigr) = \tilde
p_i ( - 1 | \tilde\sigma_t ) ,
\nonumber
\\[-8pt]
\\[-8pt]
\nonumber
&&( -1, +1) \\
\eqntext{\mbox{with probability } P_i \bigl( ( -1, +1 ) | (
\sigma_t, \tilde \sigma_t) \bigr) = p_i ( -1
| \sigma_t) - \tilde p_i ( - 1 | \tilde
\sigma_t ).}
\end{eqnarray}

Now, it is straightforward to show that under our summability condition
(\ref{B3}) on
the interaction $ J$ and $\tilde J$
and due to the boundedness of the force of the external field, the
transition probability
$P_i $ satisfies the continuity assumption (\ref{eqcontinuity}).
Hence, the decomposition
of Theorem~\ref{theodecomp} can be applied and yields the following
corollary [compare also to Theorem 3.3 of
\citet{GalGarPri10}].
%
\begin{cor}
There exists a sequence of transition probabilities $ P_k , k \geq- 1,
$ such that for any pair of symbols
$ ( a,b) \in \{ (-1,-1), (-1,1), (1,1) \} $ and any ordered pair of
configurations $ ( \sigma, \tilde\sigma) \in
\mathcal{S},$
\[
P_i \bigl((a,b) | ( \sigma, \tilde\sigma)\bigr) = \sum
_{k = -1}^\infty\lambda_i (k)
P_i^{[k]} \bigl((a,b) | ( \sigma, \tilde\sigma)
\bigl(V_i (k)\bigr) \bigr) .
\]
\end{cor}

As in Galves, L\"ocherbach and Orlandi (\citeyear{GalLocOrl10}), it can be shown that for
this decomposition,
a sufficient condition for (\ref{eqcondition2}) is
\[
\sup_{ i \in\Z^d} \sum_k \bigl| V_i (k)
\bigr| \biggl( \sum_{ j : \| j - i
\| = k } \tilde J(i,j) \biggr) < \infty
\]
and $\beta< \beta_c,$ where $\beta_c$ is solution of
\[
2 \beta\sum_{ k \geq1} \biggl( \bigl| V_i (k)\bigr |
\sum_{ j : \| j - i \| =
k } \tilde J(i,j) \biggr) = 1 .
\]
In this case, we can extend the ideas of \citet{GalGarPri10}
from the case of chains of infinite order to infinite range Gibbs measures.

Our perfect simulation algorithm simulates two ordered configurations
belonging to $\mathcal{S}$ according to the invariant distribution of
the coupled
Glauber dynamics. This yields an explicit coupling of the two Gibbs
measures $ \mu$ and~$\tilde\mu.$ Since this coupling is ordered, the
very nice argument of the proof of Theorem 3.6 in \citet{GalGarPri10} tells us
that this coupling necessarily attains the $\bar d $-distance
\[
\bar d ( \mu, \tilde\mu ) = \inf\sup_i \bigl\{ \P\bigl( \sigma(i) \neq
\tilde\sigma(i )\bigr) \dvtx (\sigma, \tilde\sigma) \mbox{ is a coupling of $\mu$ and $
\tilde\mu$} \bigr\} .
\]
Hence, our perfect simulation algorithm enables us to construct
explicit couplings achieving this distance, and as far as we know the
problem of finding explicit solutions was addressed only for finite
volume Gibbs measures up to now.

\section{Impatient user bias}

Perfect simulation procedures, very often cannot be run until the
algorithm stops, either by limitations of time or limitations of
buffer. In this section we give upper bounds for the probability of
these two types of errors.

According to our construction, the perfect simulation algorithm of
$\mu$ presented in this article is a function $ F \dvtx  [0, 1]^{ \N
\times\Z^d} $ to $S$ such that, if $(U_n)_n =( U_n (i), i \in\Z^d
)_n$ is a sequence of i.i.d. families, indexed by $\Z^d, $ of uniform
in $[0, 1] $ random variables, then for any site $i \in\Z^d, $ there
exists a stopping time $N^{(i)}_{\mathrm{STOP}}, $ such that~$F$ depends only
on the first $N^{(i)}_{\mathrm{STOP}} $ families of $(U_n)_n , $ that is, for any
measurable $B \in\mathcal{A},$
\[
P \bigl[ F \bigl( \bigl(U_1 (j)\bigr)_j, \ldots,
\bigl(U_{N^{(i)}_{\mathrm{STOP}}} (j) \bigr)_j \bigr) (i) \in B \bigr] = \mu
\bigl( \sigma(i) \in B \bigr).
\]
Note that $N^{(i)}_{\mathrm{STOP}} $ is \textit{not} the number of uniform random
variables that have to be simulated in order
to sample from $\mu; $ this number will, in general, be considerably
larger. $N^{(i)}_{\mathrm{STOP}} $ is the number of steps of the
backward sketch procedure.

A first kind of ``impatient user bias'' occurs whenever the user, for
reasons independent of the algorithm, has to stop the algorithm
after, say\vadjust{\goodbreak} $N$ steps maximal. In this case, we do not sample from $\mu
,$ but instead sample from
\[
P \bigl[ F \bigl( \bigl(U_1 (j)\bigr)_j, \ldots,
\bigl(U_{N^{(i)}_{\mathrm{STOP}}} (j) \bigr)_j \bigr) (i) \in B |
N^{(i)}_{\mathrm{STOP}} \le N \bigr] .
\]
By Proposition 6.2 of \citet{Fill} [compare also to Section 6 of
Ferrari, Fern\'andez and Garcia (\citeyear{FerFerGar02})] the error made above
can be bounded by
\[
\frac{P ( N_{\mathrm{STOP}}^{(i)} > N) }{1 -P ( N_{\mathrm{STOP}}^{(i)} > N) } \le\frac
{\gamma^N }{1 - \gamma^N }
\]
(see Theorem~\ref{theo2} above).

At each step of the backward sketch procedure, a range of order $k $ is
chosen, where $k$ is, in general, not
bounded from above. In practical situations, however, a user will be
limited in the choice of the interaction
range and will restrict the simulation to the choice of ranges bounded
by a certain upper bound $L $ that he
decided to fix in advance. More precisely, writing
\[
T_L^{(i)} := \inf\bigl\{ \tilde{T}_n^{(i)}
\dvtx \tilde K_n^{(i)} > L \bigr\} ,
\]
the user will therefore sample from the measure
\[
P \bigl[ F \bigl( \bigl(U_1 (j)\bigr)_j, \ldots,
\bigl(U_{N^{(i)}_{\mathrm{STOP}}} (j) \bigr)_j \bigr) (i) \in B | \bigl
\{N^{(i)}_{\mathrm{STOP}} \le N \bigr\} \cap\bigl\{ T_L^{(i)}
> T_{\mathrm{STOP}}^{i)} \bigr\} \bigr] .
\]
In order to control the error made induced by this ``space--time
impatient user bias,'' we have to control
\[
P\bigl( T_L^{(i)} \le T_{\mathrm{STOP}}^{i)}
\bigr).
\]
Using arguments similar to Lemma 2 of Galves, L\"ocherbach and Orlandi
(\citeyear{GalLocOrl10}), this can be bounded by
\[
P\bigl( T_L^{(i)} \le T_{\mathrm{STOP}}^{i)}
\bigr) \le\sup_{i \in\Z^d } \biggl( \frac
{ M_i - \alpha_i ( L)}{M_i} \biggr)
\frac{1}{1 - \gamma}.
\]

\section{Final comments and bibliographical discussion}

In this work we study the equilibrium measure of systems with infinite
range interactions satisfying fast decay of the long range influence
on the change rate and a certain subcriticality-criterion. For Gibbs
random fields, this regime has traditionally been studied via
cluster-expansion methods which either rely on sophisticated
combinatorial estimations [Malyshev (\citeyear{Mal80}), \citet{Sei82}, Brydges
(\citeyear{Bry86})] or inductive hypothesis and complex analysis [Koteck\'y and
Preiss (\citeyear{KotPre86}), Dobrushin (\citeyear{Dob96N1,Dob96N2})].

This is not the approach we follow here. Our approach is
probabilistic, based on an explicit construction of the dynamics and gives
probabilistic insight into the structure of the stationary law of the
process, without combinatorial or complex-analysis techniques. Let us
stress that our approach is not an alternative to cluster
expansions. It has a different regime of validity and different aims.\vadjust{\goodbreak}

Our construction is reminiscent of Harris's graphical representation
for particle systems and it is similar in spirit to procedures adopted
in Bertein and Galves (\citeyear{BerGal77}), \citet{Fer90}, \citet{vanSte99},
Ferrari, Fern\'andez and Garcia (\citeyear{FerFerGar02})
and Garcia and Mari\'c
(\citeyear{GarMar06}) among others. However, all these papers only consider
particular models, satisfying restrictive assumptions which are not
assumed in the present paper. Our approach works for any infinite
range continuous interaction under the only assumption of high-noise.

There are several techniques for perfect simulation of Markov
processes. Among the most popular ones figure \textit{coupling from the past}
(CFTP) originally proposed by \citet{ProWil96} and applied to
several special cases in a vast literature. A good review can be found
in \citet{Ken05}. This kind of technique applies to invariant
measures of Markov processes with finite coalescence time. One main
point of the CFTP technique is that one has to be able to control the
coalescence times uniformly with respect to all possible starting
points. This is an issue that becomes particularly difficult in the
case of ``big'' state spaces. The problem of large state spaces can be
overcome for processes with certain monotonicity properties or for
some specific cases. For example, for spatial point processes there is
a vast literature on the subject; we point out the works of \citet{Ken98},
Kendall and Th\"onnes (\citeyear{KenTho99}), \citet{KenMll00}
among others.

In our case, we sample directly from a time stationary realization of
the process. There is no coalescence criterion, either between coupled
realizations or between sandwiching processes. The scheme neither
requires nor takes advantage of monotonicity properties. The scheme
directly samples a finite window of the equilibrium measure in
\textit{infinite-volume}. In contrast, other CFTP algorithms [e.g.,
Kendall (\citeyear{Ken97}, \citeyear{Ken98})] focus on finite windows with fixed
boundary conditions, and the infinite-volume limit requires an
additional process of ``perfect simulation in space.'' We point out
\citet{Fer90}, \citet{van93} and \citet{vanMae94}
have also proposed construction schemes for (infinite-volume) Gibbs
measures of spin systems that can be easily transcribed into
perfect-simulation algorithms.

For continuous state spaces in systems with finite number of
components, \citet{Cai05} proposes a nonmonotone CFTP but as he points
out ``the detailed construction of the nonmonotone CFTP algorithm is
problem specific.'' Fern\'andez, Ferrari and Grynberg (\citeyear{FerFerGry07}) construct
perfect simulation for random distributions supported on a $d$
dimensional box, in particular, multivariate normal distributions
restricted to a compact set. \citet{ConKen07} show that for
a large class for positive recurrent Markov processes it is always
possible to perform CFTP, although not always feasible. In general,
for interacting particle systems with continuous state spaces, it
seems to be out of reach to apply CFTP successfully. A~recent paper by
\citet{Hub07} succeeded in using CFTP for a very specific case of a
continuous autonormal system restricted to a finite box.

The notion of random Markov chains was introduced
explicitly in \citet{Kal90} and \citet{BraKal93} and
appeared implicitly in Ferrari et al. (\citeyear{Feretal00}) and
Comets, Fern\'andez and Ferrari (\citeyear{ComFerFer02}).

\section*{Acknowledgments}

We thank Pablo Ferrari, Alexsandro Gallo, Servet Martinez and
Alexandra Schmidt for many
comments and discussions. This work is part of USP project
``Mathematics, computation, language and the brain,'' USP/COFECUB project
``Stochastic systems with interactions of variable range'' and CNPq projects
476501/2009-1 and 485999/2007-2.

%


\printaddresses

\end{document}